\documentclass[letterpaper, 10 pt, conference]{ieeeconf}  

\IEEEoverridecommandlockouts                              

\title{\LARGE \bf 
Energy-based Stabilization of Network Flows in Multi-machine Power Systems
\thanks{This work was supported in part by the ETH Z\"urich funds and the SNF Assistant Professor Energy Grant \#160573. 
}}

\author{Catalin Arghir and Florian D\"orfler%
	\thanks{C. Arghir and F. D\"orfler are with the Automatic Control Laboratory, ETH  Z\"{u}rich, 8092 Z\"{u}rich, Switzerland. Emails: 
		{\tt carghir@control.ee.ethz.ch, dorfler@ethz.ch}.}
}

\usepackage{graphicx,color}
\usepackage{animate}
\usepackage{amsmath}
\usepackage{amssymb}
\usepackage{mathrsfs}
\usepackage[dvipsnames]{xcolor}
\usepackage[vlined,ruled]{algorithm2e}

\usepackage{subfigure}

\usepackage{mathtools}
\usepackage{cite}
\usepackage{easybmat}
\usepackage{epstopdf}
\usepackage{booktabs}

\DeclareSymbolFont{bbold}{U}{bbold}{m}{n}
\DeclareSymbolFontAlphabet{\mathbbold}{bbold}
\newcommand{\vect}[1]{\mathbbold{#1}}
\newcommand{\vones}[1][]{\vect{1}_{#1}}

\newtheorem{theorem}{Theorem}[section]
\newtheorem{lemma}[theorem]{Lemma}
\newtheorem{definition}{Definition}

\newtheorem{proposition}[theorem]{Proposition}
\newtheorem{problem}{Problem}
\newtheorem{remark}{Remark}

\newtheorem{assumption}{Assumption}

\newenvironment{pfof}[1]{\vspace{1ex}\noindent{\itshape Proof of
		#1:}\hspace{0.5em}} {\hfill\QED\vspace{1ex}}

\newcommand\oprocendsymbol{\hbox{$\square$}}
\newcommand\oprocend{\relax\ifmmode\else\unskip\hfill\fi\oprocendsymbol}

\graphicspath{{./Figures/}}

\usepackage{hyperref}
 
\renewcommand{\baselinestretch}{0.99}

\begin{document}
\maketitle
\thispagestyle{empty}
\pagestyle{empty}


\begin{abstract}

This paper considers the network flow stabilization problem in power systems and adopts an output regulation viewpoint. Building upon the structure of a heterogeneous port-Hamiltonian model, we integrate network aspects and develop a systematic control design procedure. First, the passive output is selected to encode two objectives: consensus in angular velocity and constant excitation current. Second, the non-Euclidean nature of the angle variable reveals the geometry of a suitable target set, which is compact and attractive for the zero dynamics. On this set, circuit-theoretic aspects come into play, giving rise to a network potential function which relates the electrical circuit variables to the machine rotor angles. As it turns out, this energy function is convex in the edge variables, concave in the node variables and, most importantly, can be optimized via an intrinsic gradient flow, with its global minimum corresponding to angle synchronization. The third step consists of explicitly deriving the steady-state-inducing control action by further refining this sequence of control-invariant sets. Analogously to solving the so called regulator equations, we obtain an impedance-based network flow map leading to novel error coordinates and a shifted energy function. The final step amounts to decoupling the rotor current dynamics via feedback-linearziation resulting in a cascade which is used to construct an energy-based controller hierarchically.

\end{abstract}


\section{Introduction}
\label{Section: Introduction}

This paper aims to address two fundamental questions in power systems: how to explicitly characterize a steady-state operating point and how to achieve it via feedback. 
%
%
Starting from first-principles, we cast a system-theoretic treatment and arrive at the problem of stabilizing the angle configuration corresponding to an optimal network flow.

An important obstacle in stabilization of power systems dwells in the complex nonlinearity of the mechanical-to-electrical energy exchange mechanism of the rotating machines. 
Throughout the power systems literature, this aspect is mostly obviated while adopting reduced network and machine models of various degrees of fidelity \cite{c13,c18}. 

Despite the recent adoption of high order dynamical models, the analysis is still being done in a rotating frame attached to the rotor angle, a procedure difficult to extend to multi-machine networks \cite{c4,c12,c15}. 
In this work, we express these kind of models either in a stationary frame or in a rotating frame attached to a single independent angle variable for all machines. 
Thus we avoid the non-integrability obstacle encountered in \cite{c15}, while retaining the phasor interpretation, as in \cite{c10}. 

Another major impediment in control design involves the choice of equilibrium with respect to which stability of the interconnected system is referred to \cite{c12}.  
Relating to \cite{c14}, our procedure can be regarded as a reduction to a specific coupled-oscillator system, while also adopting a set-stabilization approach. 
Here, as in \cite{c18}, the dissipation and potential functions are specified based on the network topology and in terms of the machine angles.
By explicitly adressing the angle dependency in the steady-state network losses, we are able to overcome the dissipation obstacle in \cite{c18} and indeed construct a suitable energy function.

When dealing with physical systems, often complex control specifications can be formulated in terms of the level sets of an output function. 
For example, \cite{c5} provides a study of output synchronization of systems with relative degree one, as is the case here, while in \cite{c6} circular formations are targeted. 
In this spirit, we consider the consensus in angular velocity to be encoded in the passive output to be driven to zero. 
The resulting zero dynamics possess an invariant and attractive compact set on which the steady-state control action decomposes naturally into two angle-dependent components.
One component balances the steady-state dissipation throughout the network, while the other acts as a gradient of an interconnection potential energy. 
From a power systems perspective, this function encodes the canonical network objective of inductor current minimization and capacitor voltage maximization and, as we shall see, is equivalent to that of machine angle synchronization.

The rest of the paper is organized as follows. In Section \ref{Section: Prelims} notation is introduced, while in Section \ref{Section: Modelling} the power system model is presented in detail. In Section \ref{Section: problem form} we build up the control specifications from the system structure, while in Section \ref{Section: Control} we elaborate the control strategy. Before concluding, a numerical study is performed in Section \ref{Section: numerics}.

\section{Notation}
\label{Section: Prelims}

Notation will consist of the following matrices and vectors: ${j}= \begin{bsmallmatrix}0 & -1 \\ 1 & 0 \end{bsmallmatrix}$,  ${I}_2 = \begin{bsmallmatrix}1 & 0 \\ 0 & 1 \end{bsmallmatrix}$, while $\mathsf{e}_1 = \begin{bsmallmatrix} 1 \\ 0 \end{bsmallmatrix}$, $\mathsf{e}_2 = \begin{bsmallmatrix} 0\\ 1 \end{bsmallmatrix}$  are basis vectors in $\mathbb{R}^2$, ${\vones}  = \begin{bsmallmatrix} 1 & \ldots & 1\end{bsmallmatrix}^\top\in\mathbb{R}^n$, $\text{diag}(u_i)$ is a matrix with all elements of the type $u_i$ on the diagonal and zeros elsewhere. In a similar fashion we use $\text{blkdiag}(U_i)$ to lay all (square) matrices $U_i$ as blocks on the diagonal and zeros elsewhere. Denoting $\otimes$ as the Kronecker product, we make use of the symbol $\boldsymbol{j} = {I} \otimes {j}$, where ${I}$ is the identity matrix in the space of dimension implied by the context. We use the symbol $\nabla$ to denote the transpose of the Jacobian operator and $\succ$ to denote positive definiteness. If unspecified, we use index $i$ to equally refer to {each} of the $1$ up to $n$ elements, while $(\cdot, \cdot)$ is used for the column vector concatenation. Consider the angle $\theta_i \in \mathbb{S}^1$ as an object lying on the unit circle, vectorized as $\theta = (\ldots\, \theta_i, \, \ldots ) \in \mathbb{T}^{n}$, an object lying on the $n$-torus. Let $R_{\theta_i} = \begin{bsmallmatrix}\cos\theta_i & -\sin\theta_i \\ \sin\theta_i & \cos\theta_i \end{bsmallmatrix} \in \mathbb{R}^{2\times 2}$  and $\mathbf{R}_{\theta} = \text{blkdiag}(R_{\theta_i}) \in \mathbb{R}^{2n\times 2n}$ denote the corresponding rotation matrices. 
\section{MODELLING}
\label{Section: Modelling}

In this section we introduce a heterogeneous multi-machine power system model composed of $n$ synchronous machines and $m$ transmission lines which consist of inductive edges and capacitive nodes with dissipation in every element. The interconnection topology is based on a connected undirected graph defined by a signed incidence matrix $E\in\mathbb{R}^{n\times m}$. The loads of the system are modelled as constant conductances attached to the capacitive voltage buses which are associated to each synchronous machine. To preserve clarity of presentation we omit voltage busses which are not directly attached to a generator. 

\smallskip
\begin{assumption} \label{Ass:1} The three-phase system is considered to be symmetric; all resistances, inductances and capacitances have equal, positive value for each phase element.
\oprocend \end{assumption}
\smallskip

Due to this symmetry, given a three-phase quantity $\mathsf{z}_i\in\mathbb{R}^3$, we shall only consider its components restricted to the two-dimensional plane orthogonal to the vector $\vones\in \mathbb{R}^3$, known as the $\alpha\beta$-frame, as follows: ${z}_{i} = {T}_{\alpha\beta}\mathsf{z}_i$, where ${T}_{\alpha\beta}$ is given by the first two rows of the power invariant Clarke transformation, $\small{{T}_{\alpha\beta\gamma} =\begin{smallmatrix}{\sqrt{{2}/{3}}}\end{smallmatrix}\! \begin{bsmallmatrix}1 & -{1}/{2} & -{1}/{2} \\ 0 & {{\sqrt3}/{2}} & -{{\sqrt3}/{2}} \\ {1}/{\sqrt{2}} & {1}/{\sqrt{2}} & {1}/{\sqrt{2}} \end{bsmallmatrix}}$.

\subsection{Synchronous machine dynamics} 

Consider the reference the model presented in \cite{c12}, where we omit the damper windings and express the stator flux and stator current in $\alpha\beta$ coordinates. The dynamics of the $n$ synchronous machines are defined as
\begin{subequations}
\label{eq: generator-dynamics}
\begin{align}
M\dot\omega&= -D\omega - \tau_e + u_{m}  \label{eq: rotor-momenta}
\\
\dot{\theta} &= \omega	\label{eq: rotor-angle}
\\
\dot \lambda_{r}  &= -R_{r} i_{r}  + u_{r}  	\label{eq: rotor-flux}
\\
\dot \lambda_{s}  &= -R_{s} i_{s} + v  \,. 	\label{eq: stator-flux}
\end{align}
\end{subequations}
Here, $M = \text{diag}(M_{i}) \in \mathbb{R}^{n\times n}$ denotes the rotor moments of inertia, $D = \text{diag}(D_{i}) \in \mathbb{R}^{n\times n}$ denotes the rotor damping coefficients, while $\omega \in \mathbb{R}^n$, the rotor angular velocities. 
The rotor angles are represented by $\theta \in \mathbb{T}^n$, the stator winding resistances by $R_{s} = \text{diag}(R_{s_i})\otimes I_2 \in \mathbb{R}^{2n\times 2n}$ and the rotor winding resistances by $R_{r} = \text{diag}(R_{r_i}) \in \mathbb{R}^{n\times n}$. 
Each generator is connected to a voltage bus $v_i = (v_{\alpha_i}, v_{\beta_i})\in\mathbb{R}^2$ such that the stacked vector is denoted $v = (\ldots\, v_i, \, \ldots ) \in \mathbb{R}^{2n}$. 
The machines are mechanically actuated by the rotor shaft torques $u_{m}\in\mathbb{R}^n$ and the rotor winding excitation voltages $u_{r}\in\mathbb{R}^n$. Let $\lambda_{s_i} = (\lambda_{s_{\alpha_i}},\lambda_{s_{\beta_i}}) \in \mathbb{R}^2 $ represent the $\alpha\beta$ component of the stator flux in machine $i$, such that $\lambda_{s} = (\ldots\, \lambda_{s_i}, \ldots)\in\mathbb{R}^{2n}$, while $\lambda_{r} = (\ldots\, \lambda_{r_i}, \ldots)\in\mathbb{R}^{n}$ denote the rotor fluxes. 
The magnetic energy stored in all electrical machines is defined as
$$
W_{e} =  \frac{1}{2}\lambda^{\top}  \mathbf{L}_\theta^{-1}\lambda \,, \text{ where } \, \lambda = \begin{bmatrix}\lambda_s \\ \lambda_r\end{bmatrix} \, \text{ and}
$$
$$
\mathbf{L}_\theta =  \begin{bmatrix} L_s & \mathbf{R}_{\theta} (L_m  \otimes \mathsf{e}_1) \\  (L_m  \otimes \mathsf{e}_1^\top)\mathbf{R}_{\theta}^\top  & L_r\end{bmatrix} \,.
$$
Here, $L_{s} = \text{diag}(L_{s_i})\otimes I_2 \in \mathbb{R}^{2n\times 2n}$ denotes the stator self-inductances, $L_m = \text{diag}(L_{m_i}) \in \mathbb{R}^{n\times n}$, the mutual inductances and $L_r = \text{diag}(L_{r_i}) \in \mathbb{R}^{n\times n}$, the rotor self-inductances. By assuming that $L_sL_r -L_m^2\succ0$, we have that $\mathbf{L}_\theta$ is invertible and that $W_e$ is positive. We also denote the stator current by $i_{s_i} = (i_{s_{\alpha_i}}, i_{s_{\beta_i}})$, the rotor current by $i_{r_i}$ and by $\tau_{e_i}$, the electrical (air-gap) torque for each machine. We will find convenient the vectorized expressions for stator, rotor flux and electrical torque derived as
\begin{subequations}
\begin{align}
i_s = \tfrac{\partial W_{e}}{\partial \lambda_s}^\top &\Leftrightarrow \lambda_{s} =L_s i_{s} + \mathbf{R}_{\theta} (L_m  \otimes \mathsf{e}_1) i_{r} 
\label{eq: stator-currents}
\\
i_r = \tfrac{\partial W_{e}}{\partial \lambda_r}^\top &\Leftrightarrow \lambda_{r} = L_r i_{r}+ (L_m  \otimes \mathsf{e}_1^\top)\mathbf{R}_{\theta}^\top i_{s}
\label{eq: rotor-currents}
\\
\tau_{e} =  \tfrac{\partial W_{e}}{\partial \theta}^\top &\Leftrightarrow \tau_{e} = -I_{r} (L_m  \otimes \mathsf{e}_2^\top)\mathbf{R}_{\theta}^\top i_{s} \,, 
\label{eq: machine-torques}
\end{align}
\end{subequations}
where $I_r = \text{diag}(i_{r_i})$.
Finally, we will also use the vectorized expressions for stator and rotor electromotive forces
\begin{subequations}
\label{eq: machine-emf}
\begin{align}
	\xi_s &= \tfrac{\mathsf{d}}{\mathsf{d}t}\big(\mathbf{R}_{\theta} (L_m  \otimes \mathsf{e}_1) i_{r}\big)	\label{eq: s-emf}
	\\
	\xi_r &= \tfrac{\mathsf{d}}{\mathsf{d}t}\big((L_m  \otimes \mathsf{e}_1^\top)\mathbf{R}_{\theta}^\top i_{s}\big)\,,	\label{eq: r-emf}
\end{align}
\end{subequations}
representing the voltages mutually induced in the stator and rotor circuits, respectively, as per Faraday's law.

\subsection{Transmission network dynamics}

The transmission system interconnecting all generators is modelled as the following linear circuit
\begin{subequations}
\label{eq: trans}
\begin{align}
C \dot v  &= -Gv - \mathbf{E} i_{\mathsf{t}} - i_s
\\
L_{\mathsf{t}} \tfrac{\mathsf{d}}{\mathsf{d}t} i_{\mathsf{t}} &= -R_{\mathsf{t}} i_{\mathsf{t}} + \mathbf{E}^{\top} v \,,
\end{align}
\end{subequations}
where we denote by $C = \text{diag}(C_{i}) \otimes {I}_2 \in \mathbb{R}^{2n\times 2n}$ the bus voltage capacitances and by $G = \text{diag}(G_{i}) \otimes {I}_2$ the load conductances at each node with respect to the potential of the earth. Furthermore, $i_{\mathsf{t}} \in \mathbb{R}^{2m}$ denotes the vector of transmission line currents in $\alpha\beta$ coordinates, while $L_{\mathsf{t}} = \text{diag}(L_{\mathsf{t}_i}) \otimes {I}_2 \in \mathbb{R}^{2m\times 2m}$ denotes the transmission line inductances and $R_{\mathsf{t}} = \text{diag}(R_{\mathsf{t}_i}) \otimes {I}_2 \in \mathbb{R}^{2m\times 2m}$ the transmission line resistances. Here, the incidence matrix $\mathbf{E} = E\otimes {I}_2  \in \mathbb{R}^{2n\times 2m}$ describes the interconnection topology.

\begin{remark}\textbf{(The $\gamma$-component).} 
Notice that the $\alpha$ and $\beta$ components of the transmission system are independent and follow identical dynamics. The main assumption for reducing to the $\alpha\beta$-components of the three-phase circuit is that the third component, also called $\gamma$-subsystem, is not excited and follows linear, asymptotically stable dynamics.
\oprocend 
\end{remark}

\subsection{Combined system}

To write system \eqref{eq: generator-dynamics}, \eqref{eq: trans} in port-Hamiltonian form, first define the state vector as $x = (M\omega, \theta, \lambda_r, \lambda_s, Cv, L_{\mathsf{t}} i_{\mathsf{t}}) \in \mathcal{X}$ and the input vector as $u = (u_m, u_r) \in \mathcal{U}$, where the state space and input space are defined respectively as
\begin{subequations}
\label{eq: spaces}
\begin{align}
\mathcal{X} &= \mathbb{R}^n \times \mathbb{T}^n \times \mathbb{R}^{n} \times \mathbb{R}^{2n} \times \mathbb{R}^{2n} \times \mathbb{R}^{2m}
\\
\mathcal{U} &= \mathbb{R}^n \times \mathbb{R}^n \,.
\end{align}
\end{subequations}

Further consider the following Hamiltonian function
\begin{equation}
\label{eq: hamiltonian-function}
H(x) = \frac{1}{2}\omega^\top M\omega  + \frac{1}{2}\lambda ^{\top} \mathbf{L}_\theta^{-1}\lambda  + \frac{1}{2}v^{\top}Cv
 +  \frac{1}{2}i_{\mathsf{t}}^\top L_{\mathsf{t}} i_{\mathsf{t}} \,.
\end{equation}
By appropriately using constant matrices $J = -J^\top$, $R$ positive semi-definite and $B$, equations \eqref{eq: generator-dynamics}, \eqref{eq: trans} can be written in port-Hamiltonian form
\begin{subequations}
\label{eq: hamiltonian}
\begin{align}
\dot x &= (J-R) \nabla H(x) + Bu	 \label{eq: system}
\\
y &= {B}^\top \nabla H(x) \,,
\end{align}
\end{subequations}
where $\nabla H = (\omega, \tau_e, i_r, i_s, v, i_{\mathsf{t}})$ represents the co-energy variables. In this setting, we identify the passive output $y = (\omega, i_r)$ to serve as a starting point for control design.

\section{Problem Formulation}
\label{Section: problem form}

We begin with the observation that, for systems with well defined vector relative degree such as \eqref{eq: hamiltonian}, there exists a unique feedback which renders the zero-dynamics manifold invariant, see e.g. \cite{c8}. Consider $\omega_0>0$ an angular velocity to be tracked by all generators and $i_r^*\in\mathbb{R}_{>0}^n$ a vector of reference rotor currents. The associated {\it zero-dynamics manifold} can be written as
$$
\Omega = \left\{x\in\mathcal{X}: i_r = i_r^*,\, \omega = \omega_0\vones \right\} \,.
$$

The rest of this chapter develops a series of refinements of the target set $\Omega$, together with the control action which renders it invariant, arriving in the end at a suitable control specification for the power system at large.

\subsection{Zero dynamics analysis}
We will find it convenient to denote the stator impedance, bus admittance and line impedance matrices, respectively, as
\begin{equation}
\label{eq: impedances}
\begin{matrix}
	{Z}_s = R_s + \boldsymbol{j}\omega_0L_s,
	&
	{Y}_\mathsf{c} = G + \boldsymbol{j}\omega_0C, 	
	&
	{Z}_\mathsf{t} = R_\mathsf{t} + \boldsymbol{j}\omega_0L_\mathsf{t},
\end{matrix}\nonumber
\end{equation}
as well as the transmission line Laplacian, $\mathcal{L}_{\mathsf{t}} = \mathbf{E}{Z}_\mathsf{t}^{-1}\mathbf{E}^\top$. The following result, yet intermediary, will aid us in the subsequent derivations.

\begin{lemma}\textbf{(Feedback on $\Omega$).} \label{Lemma:11} 
The unique control action $u = u^\star(x)$ which renders $\Omega$ invariant for \eqref{eq: hamiltonian} is given by
\begin{subequations}
\label{eq: u-Omega}
\begin{align}
u_r^\star(x) &= R_r i_r^* + (L_m  \otimes \mathsf{e}_1^\top)\mathbf{R}_{\theta}^\top L_s^{-1}(-{Z}_s i_{s}+v)  	\label{eq: ur-omega}
\\
u_m^\star(x) &= D\omega_0{\vones} -I_r^*(L_m  \otimes \mathsf{e}_2^\top)\mathbf{R}_{\theta}^\top i_{s} \,.\label{eq: um-omega}
\end{align}
\end{subequations}
\end{lemma}
\begin{pfof}{Lemma \ref{Lemma:11}}
Using the tangency condition for a zero-level set of a function, we have that $\Omega$ is invariant {\it if and only if} the following condition holds $\forall x\in\Omega$
\begin{equation}
\label{eq: nagumo1}
	\dot\omega = 0 \,  \textit{ and } \, \tfrac{\mathsf{d}}{\mathsf{d}t} i_r = 0 \,.
\end{equation}
From \eqref{eq: rotor-currents}, \eqref{eq: r-emf} we have that $\dot\lambda_r =  \tfrac{\mathsf{d}}{\mathsf{d}t} i_r + \xi_r$ which, together with \eqref{eq: rotor-momenta}, \eqref{eq: rotor-flux} allows us to rewrite \eqref{eq: nagumo1} as
\begin{subequations}
\label{eq: nagumo2}
\begin{align}
	0 &= -D\omega - \tau_e + u_{m} \label{eq: nagumo2-m}
	\\
	\xi_r &= -R_ri_r + u_r \,. \label{eq: nagumo2-r}
\end{align}
\end{subequations}

By evaluating the expressions \eqref{eq: machine-emf} using \eqref{eq: nagumo1} we require that, $\forall x\in\Omega$,
\begin{subequations}
\label{eq: lambda-xi}
\begin{align}
	\xi_s &= \mathbf{R}_{\theta} (L_m  \otimes \mathsf{e}_2) i_{r}^*\omega_0  \label{eq: s-emf-omega}
	\\
	\xi_r &= (L_m \otimes \mathsf{e}_1^\top)\mathbf{R}_{\theta}^\top(\tfrac{\mathsf{d}}{\mathsf{d}t}i_{s}-\boldsymbol{j}\omega_0i_s)\,. \label{eq: r-emf-omega}
\end{align}
\end{subequations}
Taking a time derivative of \eqref{eq: stator-currents} yields $\dot\lambda_s =  \tfrac{\mathsf{d}}{\mathsf{d}t} i_s + \xi_s$. By using this, together with \eqref{eq: stator-flux}, in \eqref{eq: s-emf-omega}, we express \eqref{eq: r-emf-omega} as 
\begin{equation}
\label{eq: xi-r-omega}
	\xi_r =  (L_m  \otimes \mathsf{e}_1^\top)\mathbf{R}_{\theta}^\top L_s^{-1}(-{Z}_s i_{s}+v)\,,
\end{equation}
which, together with \eqref{eq: nagumo2-r} yields \eqref{eq: ur-omega}.
Finally, using $\omega\!=\!\omega_0\vones$ and $i_r \!=\! i_r^*$ in \eqref{eq: machine-torques} along with \eqref{eq: nagumo2-m}, completes the proof.
\end{pfof}

Notice that, with this feedback, the stator electro-motive force $\xi_s$ in \eqref{eq: s-emf-omega} simply becomes an Euclidean embedding of the angle $\theta$, whenever $x \in \Omega$. The zero-dynamics vector field now reads as
\begin{subequations}
\label{eq: Zero-Dynamics}
\begin{align}
\dot \theta &= \omega_0 {\vones} \label{eq: theta-dynamics}
\\
L_{s} \tfrac{\mathsf{d}}{\mathsf{d}t}i_s &= -R_s i_s + v - \xi_s
\\
C \dot v  &= -Gv - \mathbf{E}i_{\mathsf{t}} - i_s
\\
L_{\mathsf{t}} \tfrac{\mathsf{d}}{\mathsf{d}t} i_{\mathsf{t}} &= -R_{\mathsf{t}} i_{\mathsf{t}} + \mathbf{E}^{\top} v \,,
\end{align}
\end{subequations}
with $ \xi_s$ given by \eqref{eq: s-emf-omega}.

The next statement leverages the fact that, when $\xi_s$ is used to embed $\theta$ in Euclidean space, the zero-dynamics \eqref{eq: Zero-Dynamics} can be written as a block-triangular system. 
\begin{proposition}\textbf{(Stability of zero-dynamics).} \label{Prop:1} There exists a unique linear mapping $\mathbf\Pi$ such that the set $\Gamma$ defined below is globally asymptotically stable for the closed loop system \eqref{eq: hamiltonian},\eqref{eq: u-Omega}
$$\Gamma = \left\{ x\in\Omega : (i_s, v, i_{\mathsf{t}}) = \mathbf\Pi \xi_s \right\}\,.$$
\end{proposition}

\begin{pfof}{Prop \ref{Prop:1}} Notice that the closed loop dynamics \eqref{eq: hamiltonian},\eqref{eq: u-Omega}, when restricted to $\Omega$, become equivalent to \eqref{eq: nagumo1},\eqref{eq: Zero-Dynamics}. 

Using \eqref{eq: s-emf-omega} we write \eqref{eq: theta-dynamics} as $\dot\xi_s = \omega_0\boldsymbol{j}\xi_s$. By defining $z = (i_s, v, i_{\mathsf{t}})$ as the state of the three-phase electrical subsystem on $\Omega$, we can equivalently write system \eqref{eq: Zero-Dynamics} as
\begin{subequations}
\label{eq: Zero-Dynamics2}
\begin{align}
\dot \xi_s &= S \xi_s	\label{eq: exosystem}
\\
Q \dot z &= -A z + P \xi_s \,,	\label{eq: driven-system}
\end{align}
\end{subequations}
with $S = \omega_0 \boldsymbol{j} \in \mathbb{R}^{2n\times 2n}$, $Q = \text{blkdiag}(L_s, C, L_\mathsf{t})$, and
$$A = \begin{bmatrix} R_s & -{I}_{2n} & \mathsf{0} \\ {I}_{2n} & G & \mathbf{E} \\ \mathsf{0} & -\mathbf{E}^\top & R_{\mathsf{t}} \end{bmatrix}, \, P = \begin{bmatrix} -{I}_{2n} \\ \mathsf{0} \\ \mathsf{0}   \end{bmatrix}.$$

As we will see shortly, $-Q^{-1}A$ is Hurwitz, and its spectrum does not intersect the imaginary axis. Hence, there exists unique solution $\mathbf\Pi\in\mathbb{R}^{(4n+2m)\times2n}$ for the associated Sylvester equation
\begin{equation}
\label{eq: pi-map}
Q \mathbf\Pi S  =  -A \mathbf\Pi + P \,.
\end{equation}

To find $\mathbf\Pi$ we use the commutativity property of the Kronecker product, involving $\boldsymbol{j}=-\boldsymbol{j}^{-1}$. Observe that $\mathbf\Pi = (A + \omega_0 \boldsymbol{j}Q)^{-1}P$ solves \eqref{eq: pi-map} and is given by
\begin{equation}
\label{eq: pi-matrix}
\mathbf\Pi = \begin{bmatrix} - \left({Z}_s + ({Y}_\mathsf{c} + {\mathcal{L}}{_\mathsf{t}})^{-1}\right)^{-1}  \\  ({Y}_\mathsf{c} + \mathcal{L}_{\mathsf{t}})^{-1} \left({Z}_s + ({Y}_\mathsf{c} + {\mathcal{L}}{_\mathsf{t}})^{-1}\right)^{-1}  \\ {Z}_\mathsf{t}^{-1}\mathbf{E}^\top({Y}_\mathsf{c} + \mathcal{L}_{\mathsf{t}})^{-1} \left({Z}_s + ({Y}_\mathsf{c} + {\mathcal{L}}{_\mathsf{t}})^{-1}\right)^{-1}  \end{bmatrix}. 
\end{equation}

To show asymptotic stability, we define the transient component of the driven subsystem as $\tilde z = z - \mathbf\Pi\xi_s$ and, using equation \eqref{eq: pi-map}, write its dynamics as:
\begin{equation}
\label{eq: tranzient}
Q\dot{\tilde z} = -Az + P\xi_s - Q\mathbf\Pi S\xi_s = -A{\tilde z} \,.
\end{equation}

Choosing as Lyapunov function $V = \frac{1}{2}{{\tilde z}}^\top Q {\tilde z}$ yields
$$\dot V = - {{\tilde z}}^\top A {\tilde z} = -\tilde i_s^\top R_s \tilde i_s - \tilde v^\top G \tilde v -  \tilde i_{\mathsf{t}}^\top R_{\mathsf{t}} \tilde i_{\mathsf{t}}\,,$$
which is negative definite at zero, hence the transient subsystem is asymptotically stable and $-Q^{-1}A$ is Hurwitz as originally assumed. Moreover ${\tilde z}(t) \rightarrow 0$ as $t\rightarrow\infty$ here implies global asymptotic stability of $\Gamma$ for \eqref{eq: nagumo1},\eqref{eq: Zero-Dynamics}.
\end{pfof}

\begin{remark}\textbf{(Geometry of $\Gamma$).} 
We have seen that the set $\Gamma$ is positively invariant for the zero dynamics; notice that it is also compact. To illustrate its geometry, we express it as
\begin{equation}
\label{eq: geometry-gamma}
\Gamma = \left\{i_r: i_r = i_r^*\right\} \times  \left\{\omega: \omega = \omega_0{\vones}  \right\} \times \mathbb{T}^n \times \hat\pi(\mathbb{T}^n) \,, \nonumber
\end{equation}
where  we have used the {\it network flow} map $\hat\pi$, defined as 
\begin{equation}
\hat\pi : \mathbb{T}^n \rightarrow \mathbb{R}^{4n+2m} \,, \;\; \hat\pi(\theta) = \mathbf{\Pi}\mathbf{R}_{\theta} (L_m  \otimes \mathsf{e}_2) i_{r}^*\omega_0\,. \nonumber
\end{equation}
%
\oprocend 

\end{remark}

\subsection{The steady-state behaviour}
\label{subsec: steady-state control}

Consider feedback \eqref{eq: u-Omega}, suppose that system \eqref{eq: hamiltonian} is now initialized on $\Gamma$ and call such solutions $\hat x(t) \in \Gamma, \forall t>0$. Due to control invariance of $\Omega$ one has that $\omega = \omega_0\vones$ and $i_r=i_r^*$ for all time, but since we are also on the invariant set $\Gamma\subset\Omega$, we have that the solutions $\hat x(t)$ satisfy
\begin{subequations}
\label{eq: steady-state-locus}
\begin{align}
\dot{\theta} &= \omega_0\vones	\label{eq: steady-state-theta}
\\
0 &= -D\omega_0\vones - \hat \tau_e + u_{m}^\star(\hat x)	\label{eq: steady-state-omega}
\\
0 &= -R_ri_r^* + u_{r}^\star(\hat x)	\label{eq: steady-state-rotor}
\\
\boldsymbol{j}\omega_0L_s\hat i_s &=-R_s\hat i_s + \hat v - \hat \xi_s	\label{eq: steady-state-stator}
\\
 \boldsymbol{j}\omega_0C\hat v &= -G\hat v -\mathbf{E}\hat i_{\mathsf{t}} - \hat i_s	\label{eq: steady-state-voltage}
\\
\boldsymbol{j}\omega_0L_{\mathsf{t}}\hat i_{\mathsf{t}} &= -R_{\mathsf{t}}\hat i_{\mathsf{t}} + \mathbf{E}^\top \hat v \,,	\label{eq: steady-state-trans}
\end{align}
\end{subequations}
where $(\hat i_s, \hat v, \hat i_\mathsf{t})=\hat \pi(\theta)$ represent the three-phase network state variables restricted to $\Gamma$, seen as {\it phasors}, and
\begin{subequations}
\label{eq: xi-hat}
\begin{align}
	\hat\tau_e &= -I_r^*(L_m  \otimes \mathsf{e}_2^\top)\mathbf{R}_{\theta}^\top \hat i_{s}		\label{eq: tau-hat}
	\\
	\hat\xi_s &= \mathbf{R}_{\theta} (L_m  \otimes \mathsf{e}_2) I_{r}^*\omega_0\vones \,,	 \label{eq: emf-hat}
\end{align}
\end{subequations}
are the electrical torques and stator electromotive forces of system \eqref{eq: hamiltonian}, when restricted to $\Gamma$. We denote $I_r^* = \text{diag}(i_{r_i}^*)$.

Consequently, $\theta$ solely represents the state of the reduced dynamics on $\Gamma$. Notice that these dynamics can be parametrized by a constant $\theta_\mathsf{dq}\in\mathbb{T}^n$, seen as the initial condition of system \eqref{eq: steady-state-locus}, since all solutions that satisfy \eqref{eq: steady-state-locus} are of the form 
\begin{equation}
\label{eq: theta-hat}
	\theta(t) =  \omega_0\vones t + \theta_\mathsf{dq} \,.
\end{equation}

Consider the following set, which we call a {\it diagonal fiber of $\mathbb{T}^n$ through $\theta_\mathsf{dq}$},
$$\mathbb{F}|_{\theta_\mathsf{dq}}= \big\{ \theta\in\mathbb{T}^n : \theta = \omega_0\vones s + \theta_\mathsf{dq},\, s\in[0, \tfrac{2\pi}{\omega_0}) \big\} \,.$$ 
In this way, every operating point $\theta_\mathsf{dq}\in\mathbb{T}^n$ can be associated with a {\it flow} on $\Gamma$, defined as the subset
$$\Gamma|_{\theta_\mathsf{dq}} = \big\{ x\in\Omega : (i_s,v,i_\mathsf{t})\in\hat\pi(\mathbb{F}|_{\theta_\mathsf{dq}}) \big\}\,.$$
We see that such a set $\Gamma|_{\theta_\mathsf{dq}}\subset\mathcal{X}$ corresponds to a closed curve in terms of the evolution of $\hat x$, thereby characterizing a steady-state operation of our system. 

By using the $\mathbf\Pi$ map in equation \eqref{eq: xi-r-omega}, we have that $\xi_r = 0$ for the dynamics restricted to $\Gamma$. Let us now characterize how the control input \eqref{eq: u-Omega} appears on the set $\Gamma|_{\theta_\mathsf{dq}}$.

\smallskip

\begin{proposition}\textbf{(Steady-state input).} \label{Prop:3} 
Consider system \eqref{eq: hamiltonian} and a set-point $\theta_\mathsf{dq}\in\mathbb{T}^n$. The unique control action $u = u^*({\theta_\mathsf{dq}})$ which renders $\Gamma|_{\theta_\mathsf{dq}}$ invariant is given by 
\begin{subequations}
\label{eq: u-star}
\begin{align}
u_r^*({\theta_\mathsf{dq}}) &= R_r i_r^* 	\label{eq: ur-star}
\\
u_m^*({\theta_\mathsf{dq}}) &= \left(D + \mathcal{K}_\mathsf{net}({\theta_\mathsf{dq}})\right)\omega_0{\vones} \,,  \label{eq: um-star}
\end{align}
\end{subequations}
with 
$\mathcal{K}_\mathsf{net}({\theta_\mathsf{dq}}) = I_r^*(L_m\otimes \mathsf{e}_2^\top) \mathbf{R}_{\theta_\mathsf{dq}}^\top\mathbf{Y}_\mathsf{net} \mathbf{R}_{\theta_\mathsf{dq}}(L_m\otimes \mathsf{e}_2)I_r^*\,,$
and where $\mathbf{Y}_\mathsf{net} = \left({Z}_s + ({Y}_\mathsf{c} + {\mathcal{L}}{_\mathsf{t}})^{-1}\right)^{-1}$ is the equivalent admittance of the three-phase electrical subsystem.
\end{proposition}

\begin{pfof}{Proposition \ref{Prop:3}}
First, to satisfy \eqref{eq: steady-state-rotor} while on $\Gamma$, we pick $u_r^* = R_r i_r^*$. 

By letting $\hat z = (\hat i_s, \hat v, \hat i_\mathsf{t})$, denote the state of the three-phase electrical subsystem on $\Gamma$, notice that the last three equations of \eqref{eq: steady-state-locus} are equivalent to 
%
\begin{subequations}
\label{eq: pi-ss}
\begin{align}
\hat i_s &= - \left({Z}_s + ({Y}_\mathsf{c} + {\mathcal{L}}{_\mathsf{t}})^{-1}\right)^{-1} \hat \xi_s	\label{eq: is-ss}
\\
\hat z = \mathbf\Pi\hat \xi_s \, \, \, \Leftrightarrow   \, \, \,  \, \, \, \hat v &= - ({Y}_\mathsf{c} + \mathcal{L}_{\mathsf{t}})^{-1} \hat i_s	\label{eq: v-ss}
\\
\hat i_{\mathsf{t}} &= {Z}_\mathsf{t}^{-1} \mathbf{E}^\top \hat v \,. \label{eq: it-ss}
\end{align}
\end{subequations}
Using \eqref{eq: emf-hat} and \eqref{eq: is-ss} into \eqref{eq: tau-hat} allow us to express the steady-state electrical torque solely in terms of $\theta$ and $i_r^*$ as 
\begin{equation}
\label{eq: tau-e-star}
\hat \tau_e =  I_r^*(L_m\otimes \mathsf{e}_2^\top) \mathbf{R}_{\theta}^\top\mathbf{Y}_\mathsf{net} \mathbf{R}_{\theta}(L_m\otimes \mathsf{e}_2)i_r^*\omega_0  \,.
\end{equation}
Finally, if we consider the commutativity property of the Kronecker product
$$
({I}_n \otimes R_{\omega_0 t}^\top)(E \otimes {I}_2) = E \otimes R_{\omega_0 t}^\top = (E \otimes {I}_2) ({I}_m \otimes R_{\omega_0 t}^\top) 
$$
in \eqref{eq: tau-e-star}, we have that $\hat\tau_e=\mathcal{K}_\mathsf{net}({\theta_\mathsf{dq}})\omega_0{\vones}$, which is constant on $\Gamma|_{\theta_\mathsf{dq}}$. By using it in \eqref{eq: steady-state-omega}, we conclude the proof. Alternatively, this result is also obtained by using the $\mathbf{\Pi}$ map in \eqref{eq: u-Omega}.
\end{pfof}

\begin{figure}[!ht]
\centering{
\includegraphics[width=0.99\columnwidth]{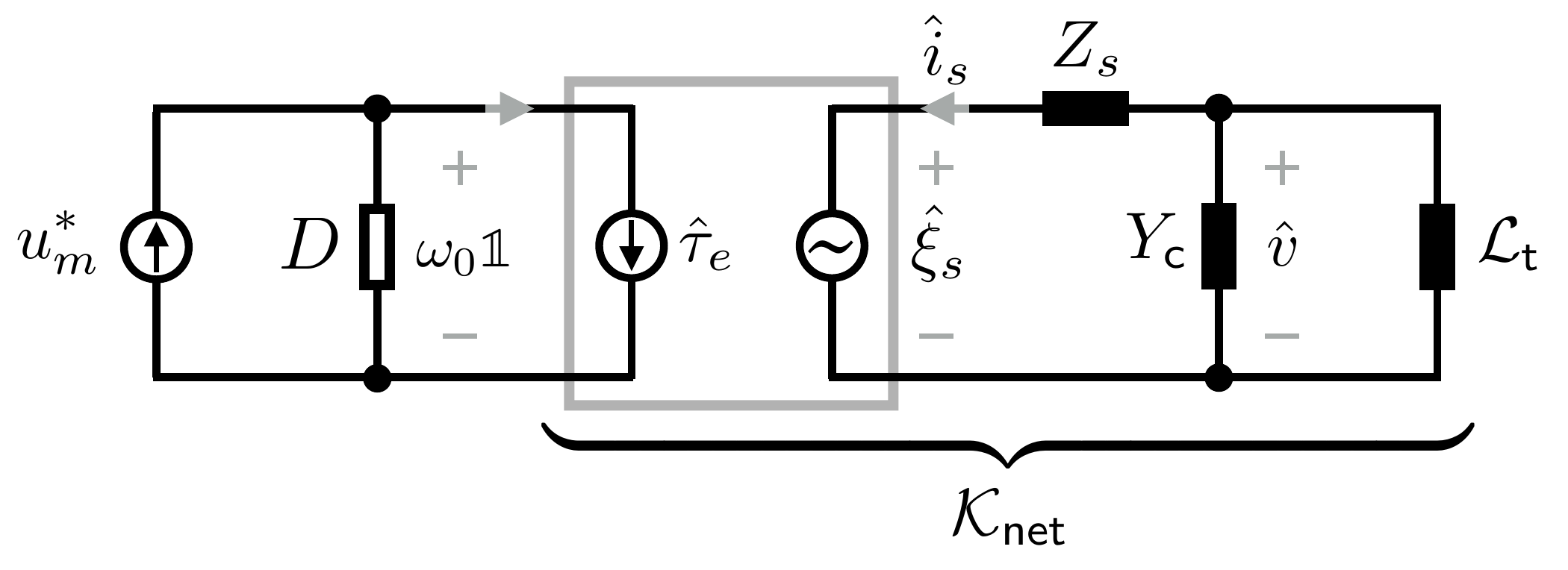}
\caption{Electrical circuit representation of the power system on $\Gamma$}
\label{Fig: ss-network}
}
\end{figure}

According to \eqref{eq: um-star}, $\mathcal{K}_\mathsf{net}$ can be interpreted as the equivalent admittance $\mathbf{Y}_\mathsf{net}$, projected via the angle $\theta_\mathsf{dq}$ as a damping factor onto the torque balance equation, as illustrated in Fig. \ref{Fig: ss-network}. We shall shortly see that this projected damping term can be decomposed into an active part and a gradient term.

\subsection{Definition of control specifications}

The aim of this section is to show how circuit aspects play a role in selecting the set-point $\theta_\mathsf{dq}$, allowing a further refinement of the control objective encoded by the set $\Gamma|_{\theta_\mathsf{dq}}$.

\begin{lemma}\textbf{(Network potential).} \label{Lemma:1} 
Consider feedback $u = u^*(\theta_\mathsf{dq})$ from \eqref{eq: u-star} and suppose that system \eqref{eq: hamiltonian} is initialized on $\Gamma|_{\theta_\mathsf{dq}}$ for some $\theta_\mathsf{dq}\in\mathbb{T}^n$. Consider the following energy function
\begin{equation}
\label{eq: S-energy}
\mathcal{S}(\theta_\mathsf{dq}) =  \frac{1}{2}\hat i_{s}^{\top} L_s \hat i_{s} - \frac{1}{2}\hat v^{\top}C \hat v + \frac{1}{2} \hat i_\mathsf{t}^{\top} L_\mathsf{t} \hat i_\mathsf{t} \,.
\end{equation}
\textit{Then}, $\theta_\mathsf{dq}$ is a critical point of $\mathcal{S}$ \textit{if and only if} $\theta_\mathsf{dq}$ is chosen such that \eqref{eq: um-star} becomes
\begin{equation}
\label{eq: input-spec2}
\resizebox{0.89\columnwidth}{!}
{$\begin{split}
u_m^* &=  I_r^{*} (L_m\otimes \mathsf{e}_2^\top) \mathbf{R}_{\theta_\mathsf{dq}}^\top\mathbf{\Pi}^\top {\begin{bsmallmatrix} R_s & & \\ & G & \\ & & R_\mathsf{t}\end{bsmallmatrix}}\mathbf{\Pi}  \mathbf{R}_{\theta_\mathsf{dq}}(L_m\otimes \mathsf{e}_2)i_r^* \omega_0 
\\&+D \omega_0\vones \,.
\end{split}$}
\end{equation}
\end{lemma}
\begin{pfof}{Lemma \ref{Lemma:1}}
One can also express $\mathbf{Y}_\mathsf{net}$ as
\begin{equation}
\label{eq: ynet}
	\mathbf{Y}_\mathsf{net}= \mathbf{\Pi}^\top \begin{bsmallmatrix}{Z}_s^\top & & \\ &{Y}_\mathsf{c} & \\ & &{Z}_\mathsf{t}^\top\end{bsmallmatrix}\mathbf{\Pi} \,.
\end{equation}
At this point, writting \eqref{eq: um-star} as
\begin{equation}
u_m^* = D\omega_0\vones + I_r^{*} (L_m\otimes \mathsf{e}_2^\top) \mathbf{R}_{\theta_\mathsf{dq}}^\top\mathbf{Y}_\mathsf{net} \mathbf{R}_{\theta_\mathsf{dq}}(L_m\otimes \mathsf{e}_1)i_r^* \omega_0 \,, \nonumber
\end{equation}
and splitting $\mathbf{Y}_\mathsf{net}$ into its symmetric and antisymmetric parts allows us to respectively decompose $u_m^*$ into an {\it active} part and a {\it gradient} term. The crucial observation is that 
one can rewrite $\mathcal{S}$ via the network flow map as
\begin{equation}
\label{eq: S-energy2}
\mathcal{S}(\theta_\mathsf{dq}) =  \frac{1}{2}\hat \xi_s^\top\mathbf{\Pi}^\top {\begin{bsmallmatrix} L_s & & \\ & -C & \\ & & L_\mathsf{t}\end{bsmallmatrix}}\mathbf{\Pi} \hat\xi_s \,,
\end{equation}
whereby taking the gradient with respect to $\theta_\mathsf{dq}$ yields
\begin{equation}
\label{eq: gradient-S}
\resizebox{\columnwidth*1}{!}{$\nabla\mathcal{S} = \omega_0I_r^{*} (L_m\otimes \mathsf{e}_2^\top) \mathbf{R}_{\theta_\mathsf{dq}}^\top\mathbf{\Pi}^\top {\begin{bsmallmatrix} L_s & & \\ & -C & \\ & & L_\mathsf{t}\end{bsmallmatrix}}\mathbf{\Pi} \mathbf{R}_{\theta_\mathsf{dq}}(L_m\otimes \mathsf{e}_1)i_r^* \omega_0 .$}\nonumber
\end{equation}
Finally, all solutions of $\nabla\mathcal{S}=0$ are critical points of \eqref{eq: S-energy}. Whenever $\theta_\mathsf{dq}$ is chosen in such a way, one has that the asymmetric part of $\mathbf{Y}_\mathsf{net}$ vanishes and \eqref{eq: um-star} becomes \eqref{eq: input-spec2}.

\end{pfof}

As we will see in the next chapter, by designing controllers to actively minimize the energy function $\mathcal{S}$, one can achieve, in terms of $(\hat i_s, \hat v, \hat i_\mathsf{t})$ variables, a saddle-point flow. In such a case, $\theta_\mathsf{dq}$ converges to a point where the energy stored in the stator impedance and in the transmission lines (edge currents) is minimized, while the energy delivered to the load (nodal voltage) is maximized. Fundamentally this is captured in the minus sign of the second term in \eqref{eq: S-energy} and can be seen as the min-flow max-cut duality in circuits. 

Furthermore, when $\theta_\mathsf{dq}$ is not a critical point of $\mathcal{S}$, the antisymmetric part of $\mathbf{Y}_\mathsf{net}$ would appear in \eqref{eq: um-star} to account for torque injections exchanged in-between generators and not delivered to the local loads. This is seen by left-multiplying \eqref{eq: um-star} by $\omega_0\vones^\top$, yielding the power balance on $\Gamma$
\begin{subequations}
\label{eq: power-balance}
\begin{align}
\omega_0\vones^\top u_m^* &=  \omega_0\vones^\top D \omega_0\vones  + \hat i_{s}^{\top} R_s \hat i_{s} + \hat v^{\top}G \hat v + \hat i_\mathsf{t}^{\top} R_\mathsf{t} \hat i_\mathsf{t} \nonumber
\\
&+ \omega_0\vones^\top \nabla\mathcal{S} \,. \nonumber
\end{align}
\end{subequations}
Since $\vones^\top \nabla\mathcal{S} = 0$ by construction, any choice of $\theta_\mathsf{dq}$ would satisfy the power balance. However, only when $\theta_\mathsf{dq} \in \nabla\mathcal{S}^{-1}(0)$, the gradient part of $u_m^*$ would be zero (element by element) in \eqref{eq: um-star}. In this case, the torque injections would exclusively account for {\it local dissipation} (local loads and the fair-share of line losses). When $\theta_\mathsf{dq}$ is chosen as a global minimizer of $\mathcal S$, then, in addition, no unnecessary circulating power would be produced in the  network, either because there are no loops in the transmission line graph or, conversely, because such a low energy level would not allow circulating currents through loops in the graph.

\begin{figure}[!h]
\centering{
\includegraphics[width=0.8\columnwidth]{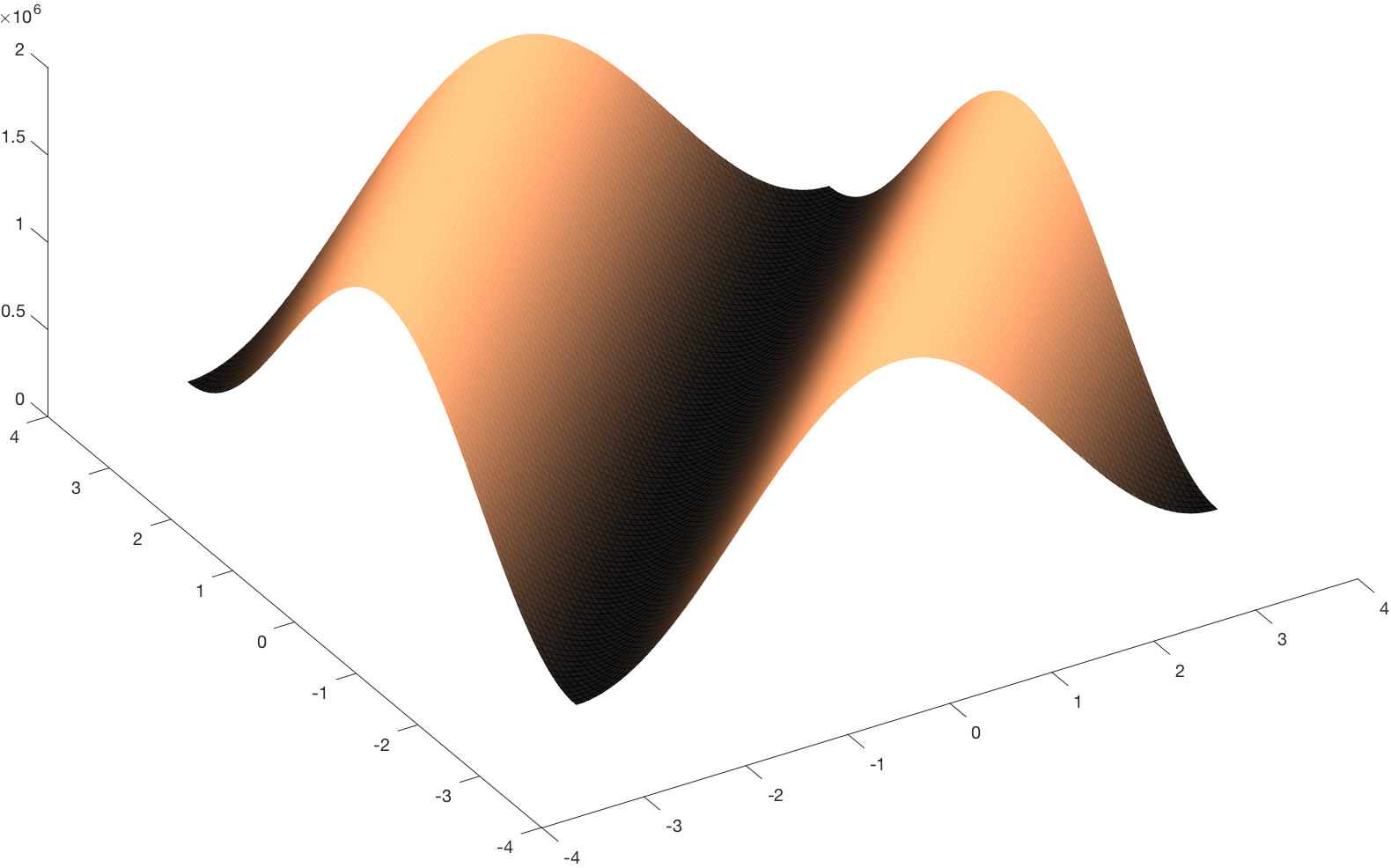}
\caption{A plot of $\mathcal{S}(\theta)$ for a two-generator test case, with  $\theta_{1},\theta_{2}\in[-\pi,\pi)$. Notice that the minimum occurs on the zero-fiber, $\mathbb{F}|_{0}\subset\mathbb{T}^2$.}
\label{Fig: gamma-torus}
}
\end{figure}

As seen in Figure \ref{Fig: gamma-torus}, for a given $\theta_\mathsf{dq}$, $\mathcal{S}(\theta)=\mathcal{S}(\theta_\mathsf{dq})$ for all $\theta\in\mathbb{F}|_{\theta_\mathsf{dq}}$. In the aim of illustrating a global stability property, we define our control specification via the set
$$\mathcal{O} = \left\{ x\in\Gamma : \nabla\mathcal{S}(\theta) = 0  \right\},$$
which is closed and bounded and, depending on the graph topology, may have more than two disconnected components. 

\section{Control Design}
\label{Section: Control}

Now we are ready to state our main control problem.

\begin{problem}\textbf{Multi-machine synchronization problem (\textit{MMSP}).} \label{Problem:2} Consider system \eqref{eq: hamiltonian} and a pair $\omega_0 > 0$, $i_r^*\in\mathbb{R}^n_{>0}$. Design 
feedback $u(x)$ such that the set $\mathcal{O}$ is globally asymptotically stable for the closed-loop system.
\end{problem}
\medskip

\begin{remark}\textbf{(Voltage control).} A general guideline for choosing the rotor current set-point $i_r^*$ is to prescribe a common stator-electromotive-force amplitude reference for all generators, denoted $e_0$, from which we derive $i_r^* = (\omega_0 L_m)^{-1}{e_0}\vones$. This baseline rotor current reference can be further adjusted to account for the voltage drop across the stator impedance and regulate instead the bus voltage amplitudes to a value closer to the set-point $e_0$. 

\oprocend 

\end{remark}

\subsection{The hierarchical control strategy}

Drawing from the distinct nature of the two inputs $u_r$ and $u_m$, 
%
we can address the system nonlinearity hierarchically,
by decoupling the rotor current and stabilizing it to a constant. 
This approach allows the integration in our framework of permanent-magnet synchronous machines (pmSM), which can be seen as having constant $i_r$, and more importantly, power electronics converters which, as shown in \cite{c9}, can be rendered equivalent to pmSM.

\begin{definition}\label{def:2}\textbf{Multi-machine synchronizing feedback (\textit{MMSF}) class.}  Consider system \eqref{eq: hamiltonian} and a pair $\omega_0 > 0$, $i_r^*\in\mathbb{R}^n_{> 0}$. Further consider inputs $u_m^*(\theta)$ from \eqref{eq: input-spec2} and $u_r^*$ from \eqref{eq: ur-star}. 
A feedback $u=(u_m, u_r)$ of the form
\begin{subequations} \label{eq: u-class}
\begin{align}
u_{r} &= u_r^* + {S}_r(x)  \label{eq: ur-class}
\\
u_{m} &= u_m^* + {S}_m(x) \label{eq: um-class}
\end{align}
\end{subequations}
is of \textit{multi-machine synchronizing feedback} (\textit{MMSF}) \textit{class if} the following conditions are satisfied:
\begin{enumerate}
\item[(i)] \emph{\textbf{Rotor current regulation:}} The set 
$$\mathcal{R} = \left\{x \in\mathcal{X}: i_r = i_r^*\right\}$$ 
is globally asymptotically stable for the closed loop;
\item[(ii)] \emph{\textbf{Velocity and angle synchronisation:}} The set 
$$\mathcal{O} = \left\{ x\in\Gamma : \nabla\mathcal{S}(\theta) = 0  \right\}$$
is globally asymptotically stable for the closed-loop dynamics restricted to the set $\mathcal{R}$.
\end{enumerate}
\oprocend \end{definition}

From the construction of feedback \eqref{eq: u-class}, we have that the stabilizing terms ${S}_r(x)$ and ${S}_m(x)$ vanish on the set $\mathcal{O}$, yielding the zero-dynamics \eqref{eq: steady-state-locus} with {associated} set of initial conditions given by the critical points of $\mathcal{S}$. 

\begin{theorem}\textbf{(Solvability of \textit{MMSP}).}\label{Thm:1} Any feedback of class \textit{MMSF} solves the \textit{MMSP} for system \eqref{eq: hamiltonian}, provided all solutions of the closed-loop system are bounded.
\end{theorem}

\begin{pfof}{Thm \ref{Thm:1}} To prove this, we set up the conditions to apply Proposition 14 in \cite{c7}, which is restated for convenience in Appendix \ref{appendix: reduction}.

\textbf{Step 1:} (Hierarchy of specifications) The sets $\mathcal{O}\subset \mathcal{R} \subset \mathcal{X}$ are closed and positively invariant for the closed-loop, which is seen from Definition \ref{def:2} and Lemma \ref{Lemma:1}. Furthermore $\mathcal{O}$ is compact since it is also contained inside the compact $\Gamma$.

\textbf{Step 2:} (Hierarchical feedback) Using Definition 4 from \cite{c7} we can see that $\mathcal{O}$ is globally asymptotically stable relative to $\mathcal{R}$, while the set $\mathcal{R}$ is in turn globally asymptotically stable for the closed-loop system.

Together with the hypothesis of bounded closed-loop solutions, we are able to use Proposition 14, case (b), in \cite{c7} and evaluate that the set $\mathcal{O}$, in this case, is globally asymptotically stable for the closed-loop system. 
\end{pfof}

\subsection{An Energy-based solution to \textit{MMSP}}

Following the strategy in Theorem \ref{Thm:1}, we construct a solution which, on one hand enforces a cascade between the rotor current dynamics and the rest of the system via feedback-linearization and on the other, uses the input designed in \eqref{eq: input-spec2} to account solely for the angle-dependent steady-state losses and induce a negative-gradient flow.

\begin{proposition}\textbf{(Solution to \textit{MMSP}).} \label{Prop:14} Consider system \eqref{eq: hamiltonian} and assume that the following holds $\forall\theta\in\mathbb{T}^n$
\begin{equation}
\label{eq: dissipation-ineq}
\resizebox{0.89\columnwidth}{!}{$D \succ \tfrac{1}{4}\omega_0^2 I_r^{*} (L_m\otimes \mathsf{e}_2^\top) \mathbf{R}_{\theta}^\top\mathbf{\Pi}^\top Q^2K^{-1}\mathbf{\Pi} \mathbf{R}_{\theta}(L_m\otimes \mathsf{e}_2)I_r^*$}
\end{equation}
where $K = \text{blkdiag}(R_s, G, R_\mathsf{t})$ and  $Q = \text{blkdiag}(L_s, C, L_\mathsf{t})$.
{\it Then}, the following controller is of \textit{MMSF} class and solves the \textit{MMSP}
\begin{subequations} \label{eq: u-energy}
\begin{align}
	u_r(x) &= R_r i_r^* + \xi_r  \,	\label{eq: ur-energy} 
	\\
	u_m(x) &= D \omega_0\vones - I_r^{*} (L_m\otimes \mathsf{e}_2^\top) \mathbf{R}_{\theta}^\top\hat i_s(\theta) - \nabla\mathcal{S}(\theta)\,.	\label{eq: um-energy}
\end{align}
\end{subequations}
\end{proposition}

\smallskip

\begin{pfof}{Proposition \ref{Prop:14}} 

\textbf{Step 1:} (Rotor current regulation) Using \eqref{eq: rotor-currents} in \eqref{eq: rotor-flux}, we can express the rotor current dynamics as:
\begin{equation}
\label{eq: rotor-current-i}
L_{r} \tfrac{\mathsf{d}}{\mathsf{d}t} i_{r} +  \tfrac{\mathsf{d}}{\mathsf{d}t}\big((L_m  \otimes \mathsf{e}_1^\top)\mathbf{R}_{\theta}^\top i_{s}\big)= -R_{r}i_{r} + u_{r} \,.
\end{equation}

By using \eqref{eq: r-emf}, we see that implementing \eqref{eq: ur-energy} assigns linear, asymptotically stable dynamics for the rotor currents
\begin{equation}
\label{eq: rotor-current-CL}
	L_{r} \tfrac{\mathsf{d}}{\mathsf{d}t} i_{r} = -R_{r}(i_{r} - i_{r}^*) \,,
\end{equation}
which implies that the set $\mathcal{R}$ is rendered globally asymptotically stable in closed-loop. 

\textbf{Step 2:} (Velocity and angle synchronization) First notice, via the decomposition of $\mathbf{Y}_\mathsf{net}$, that \eqref{eq: um-energy} is equivalent to \eqref{eq: input-spec2}. Consider now the following coordinate transformation
\begin{equation}
\label{eq: dq-transform}
\begin{smallmatrix}
\begin{split}
	\dot\theta_0&=\omega_0
	\\
	\theta_\mathsf{dq} \!&= \theta - \theta_0\vones 
	\\
	\tilde\omega \!&= \omega - \omega_0\vones 
\end{split}
& & & & &
\begin{split}
	\tilde i_{s} \!&= \mathbf{R}_{\theta_0\vones}^\top(i_s - \hat i_s(\theta)) 
	\\
	\tilde v \!&= \mathbf{R}_{\theta_0\vones}^\top(v - \hat v(\theta))
	\\
	\tilde i_\mathsf{t} \!&= \mathbf{R}_{\theta_0\vones}^\top(i_\mathsf{t} - \hat i_\mathsf{t}(\theta)) 
\end{split}
\end{smallmatrix} \,,
\end{equation}
where $\theta_0\in\mathbb{S}^1$ is an auxiliary variable and where $(\hat i_s, \hat v, \hat i_\mathsf{t})\in\mathbb{R}^{4n+2m}$ is defined as $(\hat i_s, \hat v, \hat i_\mathsf{t}) = \hat\pi(\theta)$. Secondly, assume system \eqref{eq: hamiltonian} with input \eqref{eq: u-energy} is initialized on the invariant set $\mathcal{R}$ and consider its dynamics, written in coordinates \eqref{eq: dq-transform}
\begin{equation}
\label{eq: tilda-syst}
\resizebox{1\columnwidth}{!}
{$\begin{split}
	\dot\theta_\mathsf{dq} =&~ \tilde\omega
	\\
	M\dot{\tilde\omega} =& -D\tilde\omega + I_r^*(L_m\otimes\mathbf{e}_2^\top) \mathbf{R}_{\theta_\mathsf{dq}}^\top \tilde i_{s} - \nabla\mathcal{S}({\theta_\mathsf{dq}})
	\\
	L_s \tfrac{\mathsf{d}}{\mathsf{d}t} {\tilde i}_{s} =& - {Z}_s \tilde i_{s} + \tilde v -  \mathbf{R}_{\theta_\mathsf{dq}}(L_m\otimes\mathbf{e}_2)I_r^*\tilde\omega 
	\\
	&+ \boldsymbol{j}\omega_0L_s \mathbf{Y}_\mathsf{net}\mathbf{R}_{\theta_\mathsf{dq}}(L_m\otimes\mathbf{e}_2)I_r^*\tilde\omega	
	\\
	C\dot {\tilde v} =& -{Y}_\mathsf{c} \tilde v - \tilde i_{s} - \mathbf{E}\tilde i_\mathsf{t} 
	\\
	&- \boldsymbol{j}\omega_0C({Y}_\mathsf{c} + \mathcal{L}_{\mathsf{t}})^{-1} \mathbf{Y}_\mathsf{net}\mathbf{R}_{\theta_\mathsf{dq}}(L_m\otimes\mathbf{e}_2)I_r^*\tilde\omega 
	\\
	L_\mathsf{t}\tfrac{\mathsf{d}}{\mathsf{d}t}{\tilde i}_\mathsf{t} =& -{Z}_\mathsf{t} \tilde i_\mathsf{t} + \mathbf{E}^\top \tilde v 
	\\
	&- \boldsymbol{j}\omega_0L_\mathsf{t} {Z}_\mathsf{t}^{-1}\mathbf{E}^\top({Y}_\mathsf{c} + \mathcal{L}_{\mathsf{t}})^{-1}\mathbf{Y}_\mathsf{net}\mathbf{R}_{\theta_\mathsf{dq}}(L_m\otimes\mathbf{e}_2)I_r^*\tilde\omega .
\end{split}$}\nonumber
\end{equation}
Consider now the following energy function
\begin{equation}
\label{eq: tilde-error-extended}
	\tilde{\mathcal{H}} =  \frac{1}{2}\tilde\omega^\top M\tilde\omega + \frac{1}{2}\tilde i_{s}^\top L_s \tilde i_{s} + \frac{1}{2}\tilde v ^\top C \tilde v  + \frac{1}{2}\tilde i_\mathsf{t}^\top L_\mathsf{t} \tilde i_\mathsf{t} + \mathcal{S}(\theta_\mathsf{dq}) \,, \nonumber
\end{equation}
whose derivative along the trajectories of the closed loop system reads as
\begin{equation}
\label{eq: H-dot-error}
\dot{\tilde{\mathcal{H}}} = -\begin{bmatrix}\tilde\omega \\  \tilde i_{s} \\  \tilde v  \\  \tilde i_\mathsf{t}\end{bmatrix}^\top \begin{bmatrix}D &  T_{12}(\theta_\mathsf{dq}) \\ T_{12}^\top(\theta_\mathsf{dq}) & \begin{bsmallmatrix}R_s & &  \\ & G & \\ & & R_\mathsf{t}\end{bsmallmatrix} \end{bmatrix}\begin{bmatrix}\tilde\omega \\  \tilde i_{s} \\  \tilde v  \\  \tilde i_\mathsf{t}\end{bmatrix} \,,
\end{equation}
where 
\begin{equation}
\label{eq: Q12-top}
\resizebox{1\columnwidth}{!}{$
	T_{12}^\top(\theta_\mathsf{dq}) = \tfrac{1}{2}\begin{bmatrix} -\boldsymbol{j}\omega_0L_s \mathbf{Y}_\mathsf{net} \mathbf{R}_{\theta_\mathsf{dq}}(L_m \otimes \mathbf{e}_2)I_r^* \\ \boldsymbol{j}\omega_0C({Y}_\mathsf{c} + \mathcal{L}_{\mathsf{t}})^{-1}\mathbf{Y}_\mathsf{net} \mathbf{R}_{\theta_\mathsf{dq}}(L_m \otimes \mathbf{e}_2)I_r^* \\ \boldsymbol{j}\omega_0L_\mathsf{t} {Z}_\mathsf{t}^{-1}\mathbf{E}^\top({Y}_\mathsf{c} + \mathcal{L}_{\mathsf{t}})^{-1}\mathbf{Y}_\mathsf{net} \mathbf{R}_{\theta_\mathsf{dq}}(L_m \otimes \mathbf{e}_2)I_r^*\end{bmatrix}.$} \nonumber
\end{equation}
Notice that \eqref{eq: H-dot-error} is negative definite at zero in $(\tilde\omega, \tilde i_{s}, \tilde v, \tilde i_\mathsf{t})$ {\it if and only if} the following holds $\forall\theta_\mathsf{dq}\in\mathbb{T}^n$
\begin{equation}
\label{eq: damping-condition}
D - T_{12}(\theta_\mathsf{dq})\begin{bsmallmatrix}R_s & &  \\ & G & \\ & & R_\mathsf{t}\end{bsmallmatrix}^{-1}T_{12}^\top(\theta_\mathsf{dq}) \succ 0 \,,
\end{equation}
which is equivalent to \eqref{eq: dissipation-ineq}. By LaSalle's invariance principle, we conclude that closed-loop dynamics converge to 
\begin{equation}
\label{eq: damping-set}
\mathcal{M} = \left\{(\theta_\mathsf{dq}, \tilde\omega, \tilde i_{s}, \tilde v, \tilde i_\mathsf{t}) : (\tilde\omega, \tilde i_{s}, \tilde v, \tilde i_\mathsf{t}) = 0,\, \nabla\mathcal{S}({\theta_\mathsf{dq}}) = 0 \right\}\,, \nonumber
\end{equation}
which directly implies global asymptotic stability of $\mathcal{O}$ for the closed-loop system \eqref{eq: hamiltonian}, \eqref{eq: u-energy} reduced to $\mathcal{R}$.

\textbf{Step 3:} (Boundedness) So far, we have constructed a \textit{MMSF} class controller with $S_r(x) = \xi_r$ and $S_m(x) = 0$. To see that solutions of system \eqref{eq: hamiltonian} under feedback \eqref{eq: u-energy} are bounded, first consider the rotor current dynamics \eqref{eq: rotor-current-CL} and notice that $i_r(t)$ is bounded. 

Next, consider the auxiliary variable $\psi = \mathbf{R}_{\theta} (L_m  \otimes \mathsf{e}_1) i_{r}$ which is a multiplication of two bounded terms, hence it is bounded for all time. Furthermore, consider the stator flux dynamics \eqref{eq: stator-flux} where we express, using \eqref{eq: stator-currents}, the stator current via the auxiliary variable as $ i_{s} = L_{s}^{-1}(\lambda_{s}  - \psi)$.

With this in mind we turn to the electrical subsystem, composed of \eqref{eq: stator-flux} and \eqref{eq: trans} and rewrite it as
\begin{subequations}
\label{eq: lambda-dynamics}
\begin{align}
\dot\lambda_s &= -R_sL_s^{-1}(\lambda_s -  \psi) + v 
\\
C \dot v  &= -Gv - \mathbf{E}i_{\mathsf{t}} - L_s^{-1}(\lambda_s - \psi)
\\
L_{\mathsf{t}} \tfrac{\mathsf{d}}{\mathsf{d}t} i_{\mathsf{t}} &= -R_{\mathsf{t}} i_{\mathsf{t}} + \mathbf{E}^{\top} v \,,
\end{align}
\end{subequations}
which represents a linear system, asymptotically stable in open-loop, which is driven by the bounded input $\psi$. Hence the solution of this subsystem, 
is also bounded for all time. 

We now have, using \eqref{eq: rotor-currents}, that $\lambda_{r}(t)$ is also bounded. What remains is to show that the $\omega$ component of the solution is bounded. This is seen by looking at the corresponding subsystem in closed-loop
\begin{equation}
\label{eq: omega-cl}
\begin{split}
M\dot{\omega} &= -D(\omega-\omega_0\vones) + I_r (L_m\otimes\mathbf{e}_2^\top) \mathbf{R}_{\theta}^\top  i_{s}  
\\
&+I_r^{*} (L_m\otimes \mathsf{e}_2^\top) \mathbf{R}_{\theta}^\top\mathbf{\Pi}^\top K\mathbf{\Pi} \mathbf{R}_{\theta}(L_m\otimes \mathsf{e}_2)i_r^* \omega_0 \,,
\end{split}\nonumber
\end{equation}
which is again a linear system driven by a bounded input, from which we deduce that $\omega(t)$ is bounded. Having evaluated that the entire solution $x(t)$ is bounded for all time, we can conclude that the \textit{MMSP} is solved.
\end{pfof}

Observe that the dissipation condition \eqref{eq: dissipation-ineq} resembles the one of Corollary 11 in \cite{c17}, as well as generalizes inequality (41) in \cite{c4} to the multi-machine scenario developed here. 

%
\begin{remark}\textbf{(A high-gain stabilizer)} Notice that, in controller \eqref{eq: um-energy}, the stabilizing term $S_m$ is not active. One can construct another solution, based on feedback-linearization, which fully accounts for the electrical torque extraction (when restricted to $\mathcal{R}$) and restores the negative-gradient component. In this case, the stabilizing term becomes
\begin{equation}
S_m(x) = -I_r^*(L_m  \otimes \mathsf{e}_2^\top)\mathbf{R}_{\theta}^\top(i_{s}-\hat i_s(\theta)) \,,
\end{equation}
which, compared to that of \eqref{eq: um-energy}, provides a high-gain error injection and does not require \eqref{eq: dissipation-ineq}. 
Overall, both controllers require the implementation of either the active or the gradient part of the projected dissipation term $\mathcal{K}_\mathsf{net}(\theta)$. 
\oprocend 
\end{remark}
%


\begin{remark}\textbf{(Projected gradient)} The resulting synchronising component, namely $\nabla\mathcal{S}(\theta)$ in \eqref{eq: um-energy}, has the structure of a Laplacian vector field projected onto the tangent space of $\mathbb{T}^n$. Since it acts indirectly, via the torque balance equation \eqref{eq: rotor-momenta}, the closed loop can be seen as a second-order gradient descent minimizing $\mathcal{S}$ over the angle variable $\theta$.
\oprocend 
\end{remark}

\section{Numerical experiments}
\label{Section: numerics}

To illustrate some practical aspects, we use a 3-machine model defined through the following set of parameters, all in corresponding S.I. units:  

\begin{center}
  \begin{tabular}{ | l  l  l | }
	\hline
$M_1 = 22\cdot10^3$ & $M_2 = 10^4$ & $M_3 = 45\cdot10^3$ \\ 
$D_1 = 4000$ & $D_2 = 1500$ & $D_3 = 8500$ \\ 
$L_{r,1} = 1.2$ & $L_{r,2} = 7$ & $L_{r,3} = 0.7$ \\ 
$R_{r,1} = 1.68$ & $R_{r,2} = 4.2$ & $R_{r,3} = 1.2$ \\ 
$L_{m,1} = 0.04$ & $L_{m,2} = 0.08$ & $L_{m,3} = 0.02$ \\ 
$L_{s,1} = 0.0018$ & $L_{s,2} = 0.001$ &  $L_{s,3} = 0.0066$ \\ 
$R_{s,1} = 0.166$ &  $R_{s,2} = 0.07$ &  $R_{s,3} = 0.5$ \\ 
$C_1 = 1\cdot10^{-5}$ & $C_2 = 2\cdot10^{-4}$ &  $C_3 = 4\cdot10^{-3}$ \\
$G_1 = 0.8$ & $G_2 = 0.4$ &  $G_3 = 1$ \\ 
$L_{\mathsf{t},1} = 0.0047$ &  $L_{\mathsf{t},2} = 0.0038$ &  $L_{\mathsf{t},3} = 0.0024$ \\ 
$R_{\mathsf{t},1} = 0.165$ &  $R_{\mathsf{t},2} = 0.166$ &  $R_{\mathsf{t},3} = 0.07$ \\ 
	\hline
  \end{tabular}
\end{center}

\begin{figure*}[!ht]
\centering{
\includegraphics[width=2.06\columnwidth]{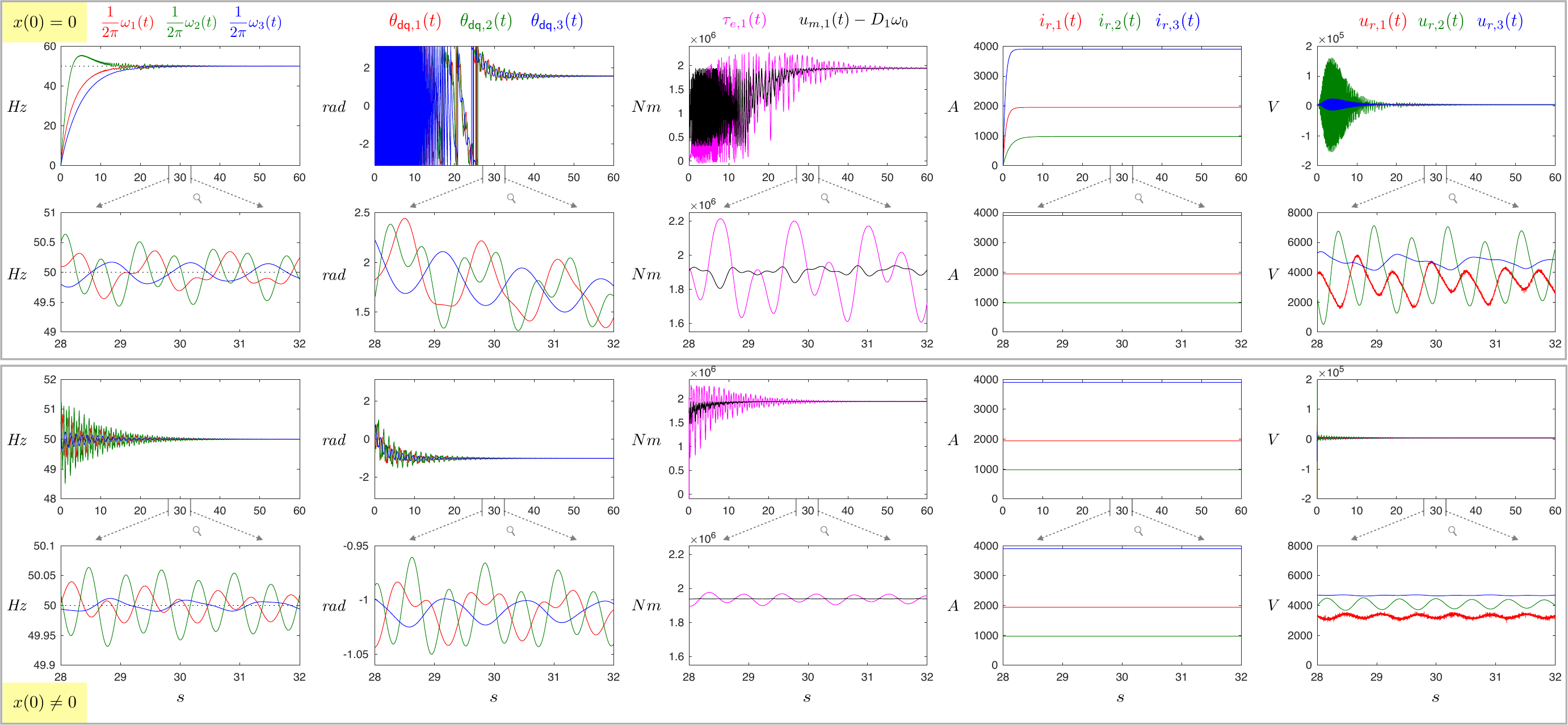}
\caption{All numbers express S.I. units, all abscissae represent time and $\theta_{\mathsf{dq},i}\in[-\pi,\pi)$. The figure shows two simulation instances: the first two rows correspond to a solution starting from zero initial condition, while the last two rows correspond to an initial condition which is closer to steady state. The variable plotted are shown on top of each column and are color coded. The second and fourth rows show an expanded view of the other two rows.}
\label{Fig: time-plots}
}
\end{figure*}

The parameters of the first synchronous generator were inspired from Example 3.1 in \cite{c20}, while the other two machines are slight variations of it. Furthermore, the topology encoded by $E=\begin{bsmallmatrix}-1&1&0 \\ 0&-1&1 \\ 1&0&-1\end{bsmallmatrix}$ is a complete graph with $n\!=\!m\!=\!3$. The system is initialized either as $x(0)=0$, or as $\omega(0)\!=\! (0.99,1.01,0.999)\,\omega_0$,  $\theta(0) \!=\! (0,-\pi/4,\pi/4)$, $i_r(0)\!=\!i_r^*$, $i_s(0)\!=\!(1.47,1.02,1.53,1.04,0.57,0.37)\cdot10^4$, $v(0)\!=\!(-1.64,-1.06,-1.53,-1.51,-1.79,-0.69)\cdot10^4$ and $i_\mathsf{t}(0)\!=\!(-0.29,-0.11,0.65,0.31,-0.48,-0.24)\cdot10^4$. The output reference was chosen as $\omega_0 \!=\! 2\pi\cdot50$, $i_{r,1}^* \!=\! 1950$, $i_{r,2}^* \!=\! 975$, $i_{r,3}^* \!=\! 3900$. Observe in Figure \ref{Fig: time-plots}, that the output regulation specifications are met, along with boundedness of solutions. In addition, the angles synchronize to zero relative differences: this is due to the fact that, for a typical choice of parameters, the syncrhonizing potential $\mathcal{S}$ has the global minimum corresponding to identical angles.

Notice that the angular velocity overshoot of generator $2$ is induced by the long transient of its rotor current. This has been created with the purpose of showing the effects of hierarchical control. 
From an implementation perspective it is worth noting that the small oscillations of the control action, during mild transients, would have to be tracked by a sufficiently fast torque actuation.

\section{Conclusions}
\label{Section: Conclusions}

In this article, starting from a first-principle power system model, we develop a series of invariant sets in the aim of characterizing its steady-state behavior. 
On these sets, the steady-state control action is subsequently derived, resulting in a gradient vector field acting on the machine torque balance equation. 
Through the use of a steady-state network map, parametrized in terms of the angle variable, we decompose the equivalent admittance of the three-phase circuit to uncover a canonical potential function. 
This function incorporates the effects of bus admittance and stator and line impedances, thus representing an extension to that of the classical swing-equations.
By accounting specifically for the steady-state losses, we arrive at a coordinate transformation and an associated shifted Hamiltonian which allow us to design an energy-based controller.
As feedback, the controller requires machine angle measurements as well as full knowledge of network parameters.
For this synthesis, we propose a framework which hierarchically stabilizes the excitation current first, and then simultaneously achieves frequency consensus and an angle configuration corresponding to the optimizers of the a network potential function. 

The discussion on robustness of implementation and how to achieve a non-identical angle configuration, corresponding to generators exchanging power over the transmission lines, is deferred to a future publication.

\renewcommand{\baselinestretch}{0.99}
\bibliographystyle{IEEEtran}
\bibliography{alias,Main,FB,new}

\appendix

\subsection{Elements of reduction theory}
\label{appendix: reduction}

We define the dynamical system $\dot x = f(x)$ and $\phi(t, x_0)$, its solution at time $t$ with initial condition $x_0$. Furthermore $B_\epsilon(x_0)$ denotes the ball or radius $\epsilon$ centred at $x_0$. 
\begin{definition}\textbf{(Definition 4 from \cite{c7}).} Let $\Gamma_1$ and $\Gamma_2$, $\Gamma_1\subset\Gamma_2\subset\mathcal{X}$, be closed positively invariant sets. We say that $\Gamma_1$ is \textit{stable relative to} $\Gamma_2$ if, for any $\epsilon>0$, there exists a neighbourhood $\mathcal{N}(\Gamma_1)$ such that $\phi(\mathbb{R}^+, \mathcal{N}(\Gamma_1) \cap \Gamma_2) \subset B_\epsilon(\Gamma_1)$. Similarly, one modifies the notion of (global) asymptotic stability by restricting initial conditions to lie on $\Gamma_2$.
\oprocend \end{definition}

\smallskip
\begin{proposition}\textbf{(Proposition 14, case (b), from \cite{c7}).} \label{Prop:5} Consider system $\dot x = f(x,u)$, and assume that there exists a locally Lipschitz feedback $\bar u(x)$ making the sets $\Gamma_1 \subset \ldots \subset \Gamma_l$ positively invariant for the closed-loop system. Let $\Gamma_{l+1} = \mathcal{X}$, and consider the following conditions for the closed-loop system $\dot x = f(x,\bar u(x))$:
\begin{enumerate}
\item[(i')] For $i = 1,\ldots, l$, $\Gamma_i$ is globally asymptotically stable relative to $\Gamma_{i+1}$ for the closed-loop system.
\item[(iii)] All trajectories of the closed-loop system are bounded.
\end{enumerate}
Then, the following implication holds:
\begin{enumerate}
\item[(b)] $(\text{i'}) \land (\text{iii}) \land (\Gamma_1 \text{ is compact}) \Rightarrow \Gamma_1$ is globally asymptotically stable for the closed-loop system.
\end{enumerate}
\end{proposition}


\end{document}